\documentclass[journal]{IEEEtran}

\usepackage{setspace}
\usepackage[english]{babel}
\usepackage{amsmath}
\usepackage{mathalfa}
\usepackage{mathrsfs}

\usepackage{amssymb}
\usepackage{amsthm}
\usepackage{color}
\newtheorem{asm}{Assumption}

\newtheorem{thm}{Theorem}

\usepackage{graphicx}
\graphicspath{ {pics/} }
\allowdisplaybreaks

\IEEEoverridecommandlockouts  

\begin{document}
	
	\title{Extremum Seeking for Traffic Congestion Control with a Downstream Bottleneck}
	
	\author{Huan~Yu,~\IEEEmembership{Member,~IEEE,} Shumon Koga,~\IEEEmembership{Member,~IEEE,} Tiago Roux Oliveira,~\IEEEmembership{Senior Member,~IEEE,} and\\ Miroslav Krstic,~\IEEEmembership{Fellow,~IEEE}%

	\thanks{Huan Yu, Shumon Koga and Miroslav Krstic are with the Department of Mechanical and Aerospace Engineering,
	University of California, San Diego, 9500 Gilman Dr, La Jolla, CA 92093, United States. {(email: huy015@ucsd.edu; skoga@ucsd.edu; krstic@ucsd.edu)}}%
\thanks{Tiago Roux Oliveira is with the Department of Electronics and Telecommunication Engineering, State University of Rio de Janeiro (UERJ), Rio de Janeiro 20550-900, Brazil.{ (email: tiagoroux@uerj.br)}}

	\thanks{The third author would like to thank the Brazilian funding agencies CNPq, CAPES, and FAPERJ for the financial support.}}

\maketitle

\begin{abstract}
	This paper develops boundary control for freeway traffic with a downstream bottleneck. Traffic on a freeway segment with capacity drop at outlet of the segment is a common phenomenon leading to traffic bottleneck problem. The capacity drop can be caused by lane-drop, hills, tunnel, bridge or curvature on the road. If incoming traffic flow remains unchanged, traffic congestion forms upstream of the bottleneck due to outgoing traffic overflowing its capacity. Therefore, it is important for us to regulate the incoming traffic flow of the segment so that the outgoing traffic at the bottleneck can be discharged with the maximum flow rate. Traffic densities on the freeway segment are described with Lighthill-Whitham-Richards (LWR) macroscopic Partial Differential Equation (PDE) model. To prevent the traffic congestion forming upstream of the bottleneck, incoming flow at the inlet of the freeway segment is controlled so that the optimal density could be achieved to maximize the outgoing flow and not to surpass the capacity at outlet. The density and traffic flow relation, described with fundamental diagram, is assumed to be unknown at the bottleneck area. We tackle this problem using Extremum Seeking (ES) Control with delay compensation for LWR PDE. ES control, a non-model based approach for real-time optimization, is adopted to find the optimal density for the unknown fundamental diagram. A predictor feedback control design is proposed to compensate the delay effect of traffic dynamics in the freeway segment. In the end, simulation results validate a desired performance of the controller on the nonlinear LWR model with an unknown fundamental diagram.	
\end{abstract}


\IEEEpeerreviewmaketitle
	
\section{Introduction}

    When there are uphills, curvature or lane-drop further downstream on freeway, a bottleneck with lower capacity could appear. Traffic congestion then forms upstream of the bottleneck if there is no traffic regulation. Ramp metering and variable speed limit (VSL) have been proved to be very effective in freeway traffic management system. Boundary control of traffic with ramp metering or VSL in presence of downstream bottleneck is studied in this work.

  The first local ramp metering strategy that was proposed based on feedback control theory is ALINEA developed by \cite{Papa1}, and later on an adaptive strategy was employed by AD-ALINEA when downstream occupancy is uncertain~\cite{Papa}. The traffic flow entering the freeway is controlled from ramp metering on-ramps so that the downstream mainline traffic flow is maximized locally or the optimal freeway network traffic is achieved coordinately. The ALINEA algorithm uses real-time measurement of downstream occupancy and its set point value to calculate their difference, and then the control input is designed via the integration of the errors over time. In the presence of a downstream bottleneck, Proportional-Integrator (PI) ALINEA was developed by~\cite{Wang} to improve performance of the closed-loop system. A comparative study by \cite{Kan} is conducted in comparison with ALINEA. PI-ALINEA is proposed as an extension to ALINEA by measuring the downstream bottleneck occupancy and feeding it back to the local ramp-metering. In \cite{Wang}, the stability of the closed-loop system with PI-ALINEA, a discretized ordinary differential equation (ODE) system, is proved with Lyapunov analysis. Simulation demonstrated that PI-ALINEA improved significantly than ALINEA in the case of distant downstream bottleneck. 

Control of lane-drop bottleneck by VSL was explored by~\cite{Jin}. Authors approximated LWR model, which is a first-order macroscopic PDE model, with the discretized ODE link queue model. A Proportional-Integrator-Derivative (PID) controller is employed for VSL control strategy. 

The control of traffic with lane-drop problem modeled with macroscopic LWR PDE was firstly investigated by \cite{BL2}. The traffic dynamics on a stretch of freeway upstream of the bottleneck area is governed by LWR model. The predictor feedback control law is designed for the ramp metering at the inlet of the freeway so that the density at bottleneck area is regulated to a desired equilibrium.  
This work assumes the prior knowledge of the optimal density that could maximize the discharging flow at the bottleneck area. However, the density and traffic flow relation at bottleneck area is usually hard to obtain or estimate, especially when the bottleneck is caused by a random accident and the traffic needs to be regulated immediately. 

In this paper, we consider a freeway segment with bottleneck located at the outlet where road capacity drops. The traffic dynamics of the freeway segment is described with LWR model. The density-flow relationship at the bottleneck area is described with a nonlinear map at the outlet where the optimal density is unknown. We apply ES control, a non-model based real-time adaptive control technique, in order to find the unknown optimal density at the bottleneck. Since the control is actuated from the upstream freeway of the bottleneck, the delay effect of the traffic dynamics is compensated in designing ES control.

ES control has been intensively studied over the recent years ~\cite{Ariyur,benosman,Jan,Gha, Gha1,Gha2,Guay1,MK6,MK5,MK2,Tiago,Paz, Rusiti, Tan,Tan2,Wang}, especially after the theoretical work by \cite{MK5} proving the convergence of cost function to a neighborhood of the optimal value by means of averaging analysis and singular perturbation. 
ES approach relies on a small periodic excitation, usually sinusoidal to disturb the parameters being tuned and the effect of the parameters is then quantified by the output of a nonlinear map. The search of the optimal value is therefore generated. Despite the large number of previous work on ES control, authors in \cite{Tiago} firstly considered the problem of ES control in the presence of delays. The proposed method is based on the predictor-based feedback for delay compensation of \cite{Tiago}. Using backstepping transformation by \cite{MK4} and averaging in infinite dimensional systems in \cite{Hale}, the stability analysis is rigorously obtained. The averaging based approach is employed due to the need to estimate the unknown second-order derivative of nonlinear map at bottleneck. 
  
Our contribution lies in the following aspects: this is the first work on control of traffic governed by LWR PDE model in the presence of unknown downstream bottleneck. The optimal density input at inlet of the freeway segment is achieved by estimating the unknown nonlinear map at the outlet. ES control with delay compensation is firstly adapted to this traffic problem. The traffic dynamics is represented with linearized LWR model in the theoretical analysis, but the simulation is conducted on the nonlinear LWR model and ES control design is validated for the nonlinear system.

The outline of this paper: 	
we firstly introduce LWR PDE model for the freeway segment upstream of bottleneck and describe density-flow relation at bottleneck with a nonlinear map. For the linearized error system, we design a predictor feedback control law with delay compensation. Stability analysis is conducted for the closed-loop system using backstepping transformation and averaging approach. To illustrate our result, simulation is performed on the nonlinear LWR PDE model and a quadratic fundamental diagram is considered. The conclusion and discussion of future work are given in the end.

	
\section{Problem Statement}	
We consider a traffic problem on a freeway-segment with lane drop downstream of the segment. The freeway segment upstream of the bottleneck and the lane-drop area are shown in Fig.1 which illustrates the clear Zone C and the bottleneck Zone B respectively. To prevent the traffic in Zone B overflowing its capacity and then causing congestion in the freeway segment, we aim to find out the optimal density ahead of Zone C that maximizes outgoing flux of Zone B given unknown density-flow relation. Traffic dynamics in Zone C is described with macroscopic traffic model for aggregated values of traffic density. The traffic dynamics in lane-drop Zone B is usually difficult to describe with mathematical model and thus assumes that the fundamental diagram is unknown. 

Here we choose the first-order LWR model instead of more sophisticated second-order model, e.g. Aw-Rascle-Zhang (ARZ) by~\cite{AW, Zhang}, for the following reasons. The second-order ARZ model consisting of two PDEs governs both the traffic density and velocity. The second-order ARZ PDE model was brought up to resolve the issue that equilibrium fundamental diagram are not suited for the congested regime. Traffic velocity could vary from the single-valued equilibrium function. In this problem, we consider the free regime for freeway segment to prevent the formation of traffic congestion in bottleneck area. LWR model therefore suit our needs.

 We control traffic flow entering at the inlet of Zone C upstream of block Zone B. Traffic dynamics of Zone C is described with the first-order LWR model. Therefore, we design ES control for an unknown static map with actuation dynamics governed by a nonlinear hyperbolic PDE. The control objective is to find the optimal input density at inlet of Zone C that drives the measurable output flux of Zone B to its unknown optimal value of an unknown fundamental diagram.

	\subsection{LWR Traffic Model} 	
	The traffic dynamics in Zone C upstream of Zone B is described with the first-order, hyperbolic LWR model. 
	\begin{figure}[t!]
		\centering
		\includegraphics[width=0.35\textwidth]{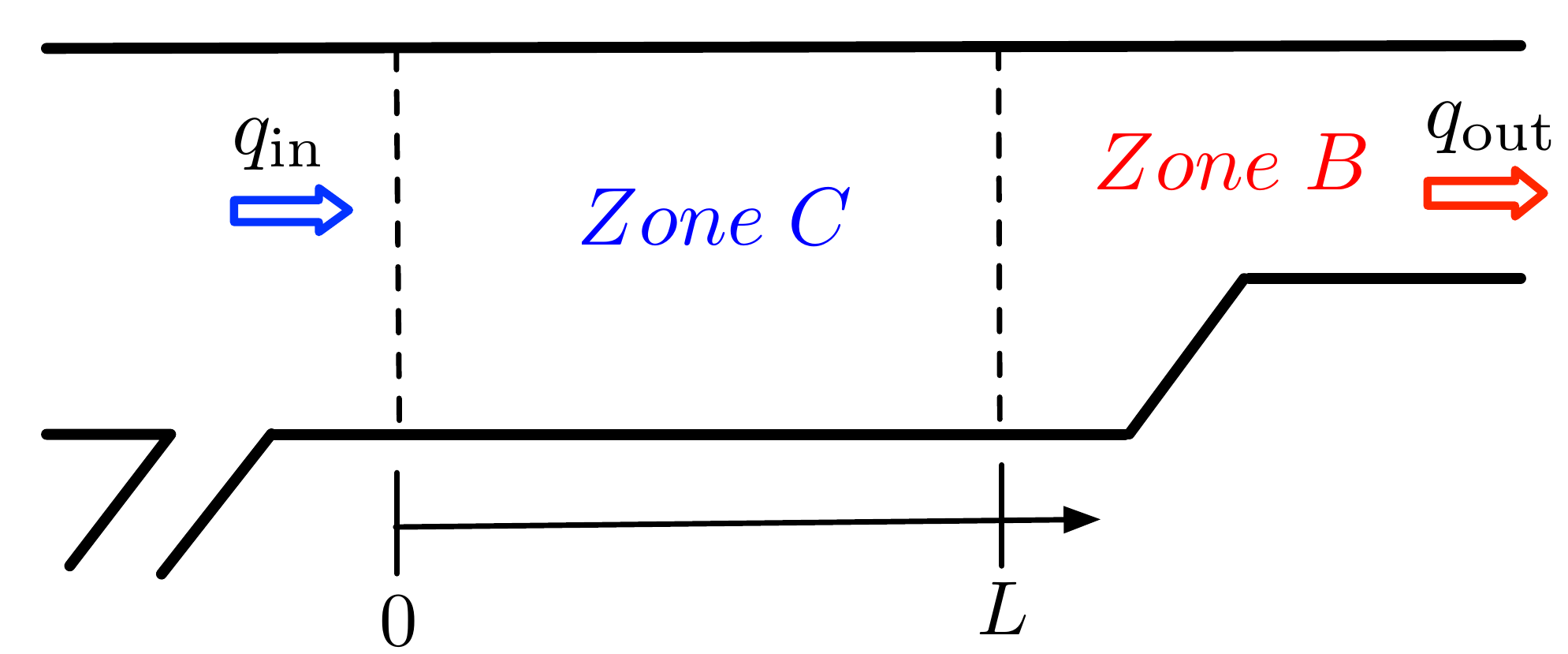}
		\caption{Traffic on a freeway segment with lane-drop}
	\end{figure}
	 Traffic density $\rho(x,t)$ in Zone C is governed by the following nonlinear hyperbolic PDE, where $x\in[0,L]$, $t\in[0,\infty)$, 
\begin{align}\label{sys1} 
\partial_t \rho + \partial_x( \rho V(\rho))=&0, 
\end{align}
where traffic velocity follows an equilibrium velocity-density relation $V(\rho)$. The fundamental diagram of traffic flow and density function $Q(\rho)$ is given by 
\begin{align}
	Q(\rho) = \rho V(\rho).
\end{align}
There are different models to describe the flux and density relation. A basic and popular choice is Greenshield's model for $V(\rho)$ which is given by 
\begin{align}
V(\rho)&=v_f\left(1-\frac{\rho}{\rho_m}\right), \label{vf}
\end{align}
where $v_f \in \mathbb{R^+}$ is defined as maximum velocity and $\rho_m \in \mathbb{R^+}$ is maximum density for Zone C. Then the fundamental diagram of flow and density function $Q(\rho)$ is in a quadratic form of density,
\begin{align}
Q(\rho) = -\frac{v_f}{\rho_m}\rho^2 + v_f \rho.
\end{align}
In practice, quadratic fundamental diagram does not provide a good fitting with real traffic density-flow data. The critical density usually happens at $20\% $ of the maximum value of density. There are several other equilibrium model e.g. Greenberg model, Underwood model and diffusion model for which fundamental diagrams are nonlinear functions. According to Taylor expansion, second-order differentiable nonlinear function can be approximated as a quadratic function in the neighborhood of its extremum. The following assumption is made for the nonlinear fundamental diagram. The stability results derived in this paper holds locally for the general form of $Q(\rho)$ that satisfy the following assumption. 

\begin{asm} 
	We assume the fundamental diagram $Q(\rho)$ is a $C^2$ function, then $Q(\rho)$ can be decomposed at the critical density $\rho_{c}$ as follows:
	\begin{align}
		Q(\rho) = q_c + \frac{Q^{\prime\prime}(\rho)}{2}(\rho(t)-\rho_c)^2,
	\end{align}
	where $q_c = Q(\rho_c)$ defined as the road capacity or maximum flow, with assumption that $Q^{\prime\prime}(\rho)< 0$ is satisfied.
\end{asm}

\subsection{Downstream Bottleneck Problem}
When there is a bottleneck present downstream, the density at outlet of Zone C is $\rho(L,t)$ governed by PDE in \eqref{sys1} for $x\in[0,L]$, $t\in[0,\infty)$ and boundary condition at inlet in \eqref{bc}. 
The inlet boundary flow is, 
\begin{align} \label{bc}
q_{\textrm{in}}(t) = Q(\rho(0,t)).
\end{align} 
The control objective is to design the traffic flow input $q_{\textrm{in}}(t)$ so that the outgoing flow in lane-drop area Zone B $q_{\textrm{out}}(t)$ is maximized.

The traffic dynamics in Zone B is described by an unknown fundamental diagram since density and traffic flow relation at the bottleneck area is hard to determine. Therefore, we assume that the equilibrium fundamental diagram for Zone B is an unknown quadratic map $Q_B(\rho)$, shown in Fig.2.	
The measurement of traffic flow in Zone B, $q_{\textrm{out}}(t)$ is defined by $Q_B(\rho)$ with outlet density $\rho(L,t)$ at outlet,
	\begin{align}
		q_{\textrm{out}}(t) =& Q_B(\rho(L,t) ).
	\end{align}
Due to the lane-drop at outlet, maximum density and road capacity reduced at Zone B compared with Zone C. We consider that optimal density of Zone B is smaller than critical density $\rho^{\star}\in \mathbb{R^+}$ of Zone C.
We aim to find out the critical outlet density $\rho^{\star}$ of Zone C that maximize $q_{\textrm{out}}(t)$ in Zone B,
\begin{align}
	q_{\textrm{out}}(t) = q^\star + \frac{H}{2}(\rho(L,t)-\rho^\star)^2, \label{static}
\end{align}	
where $q^\star \in \mathbb{R^+}$ is the unknown optimal output flow for Zone B and $H<0$ is the unknown Hessian of the static map $Q_B$.
\begin{figure}[t!]
	\centering
	\includegraphics[width=0.38\textwidth]{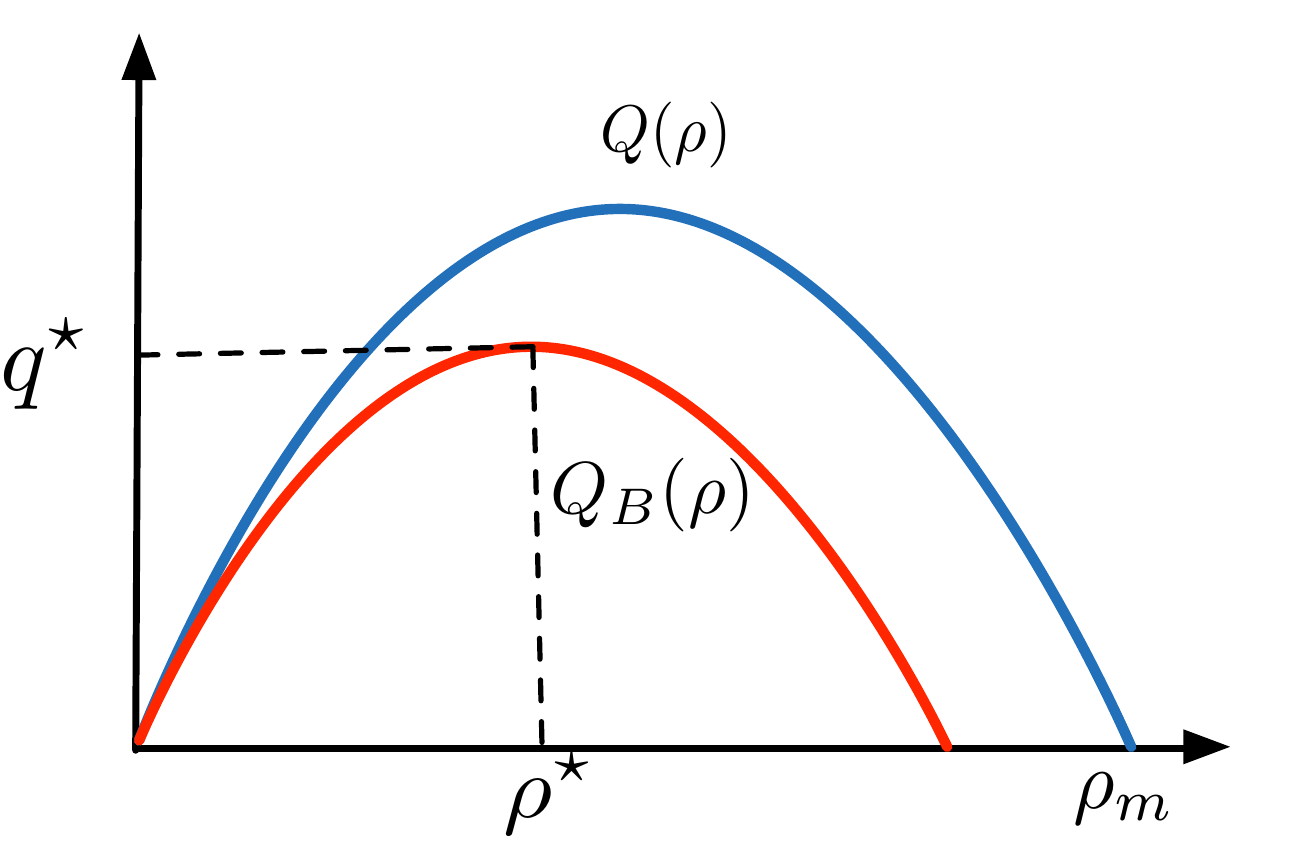}
	\caption{Quadratic fundamental diagram with distant bottleneck}
\end{figure}
\subsection{Linearized Reference Error System} 
We linearize the nonlinear LWR model around a constant reference density $\rho_{r} \in \mathbb{R^+}$, which is assumed to be close to the optimal density $\rho^{\star}$.
Note that the reference density $\rho_{r}$ is in the free regime of $Q(\rho)$ of Zone C thus is smaller than the critical density $\rho_{c}$ and therefore the following is satisfied
 \begin{align}
 	\rho_{r} <\rho_c.
 \end{align}
  Here we do not specify fundamental diagram $Q(\rho)$ for Zone C but require assumption 1 to be satisfied.
Define the reference error density as
\begin{align}
	\tilde{\rho}(x,t) = \rho(x,t) - \rho_{r},
\end{align}
and reference flux $q_r$ is 
\begin{align}
	q_{r} =& Q(\rho_r) > 0.
\end{align}
By the governing equation \eqref{sys1} together with \eqref{vf}, the linearized reference error model is derived as 
\begin{align}
\partial_{t} \tilde{\rho}(x,t)   + u \partial_{x} \tilde{\rho}(x,t)  =& 0,  \label{trho}\\
\tilde{\rho}(0,t) =& \rho(0,t)- \rho_{r}, \label{in}
\end{align} 
where the constant transport speed $u$ is given by
\begin{align} 
\notag u =& Q^\prime ( \rho)|_{\rho =\rho_r } \\
= &V(\rho_r) + \rho_r V^\prime ( \rho)|_{\rho =\rho_r }.
\end{align} 
The equilibrium velocity-density relation $V(\rho)$ is a strictly decreasing function. 
The reference density $\rho_{r}$ is in the left-half plane of the fundamental diagram $Q_c(\rho)$ which yields the following inequality for the propagation speed $u$,
\begin{align}
	u>0.
\end{align}
According to \eqref{bc} and \eqref{in}, we define the input density as
\begin{align}
\varrho(t)  = &{\rho}(0,t),
\end{align}
and the linearized input at inlet is 
\begin{align}
\tilde \varrho(t)  =& \varrho(t) - \rho_{r}.
\end{align}
The linearized error dynamics in \eqref{trho}, \eqref{in} is a transport PDE with an explicit solution for $t > \frac{x}{u}$ and thus is represented with input density
\begin{align} 
\tilde{\rho}(x,t) = \tilde\varrho\left( t - \frac{x}{u} \right),
\end{align} 
The density variation at outlet is
\begin{align}
	\tilde{\rho}(L,t) = &\tilde \varrho\left( t - D \right).\label{orho1}
\end{align}
where the time delay $D$ is defined as 
\begin{align} 
D = \frac{L}{u}.
\end{align} 
Therefore, the density at outlet is given by a delayed input density variation and the reference
\begin{align} 
{\rho}(L,t) = & \rho_r + \tilde {\rho}(L,t).\label{orho2}
\end{align} 
Finally, substituting \eqref{orho1}, \eqref{orho2} into the static map \eqref{static}, we arrive at the following
\begin{align}
	q_{\textrm{out}}(t) =&  q^\star + \frac{H}{2}( \tilde\varrho\left( t - D \right) + \rho_{r} -\rho^{\star} )^2 \notag\\
	= &  q^\star + \frac{H}{2} \left( \varrho\left( t - D \right)  - \rho^{\star} \right)^2.  \label{out}
\end{align}	
The control objective is to regulate the input $q_{\textrm{in}}(t)$ so that $\varrho\left( t - D \right)$ reaches to an unknown optimal $\rho^\star$ and the maximum of the uncertain quadratic flux-density map $q_{\textrm{out}}(t)$ can be achieved. We can apply the method of extremum seeking for static map with delays developed in~\cite{Tiago}. The extremum seeking control is designed for finding the extremum of the unknown map. 

\begin{figure}[t!]
	\centering
	\includegraphics[width=0.30\textwidth]{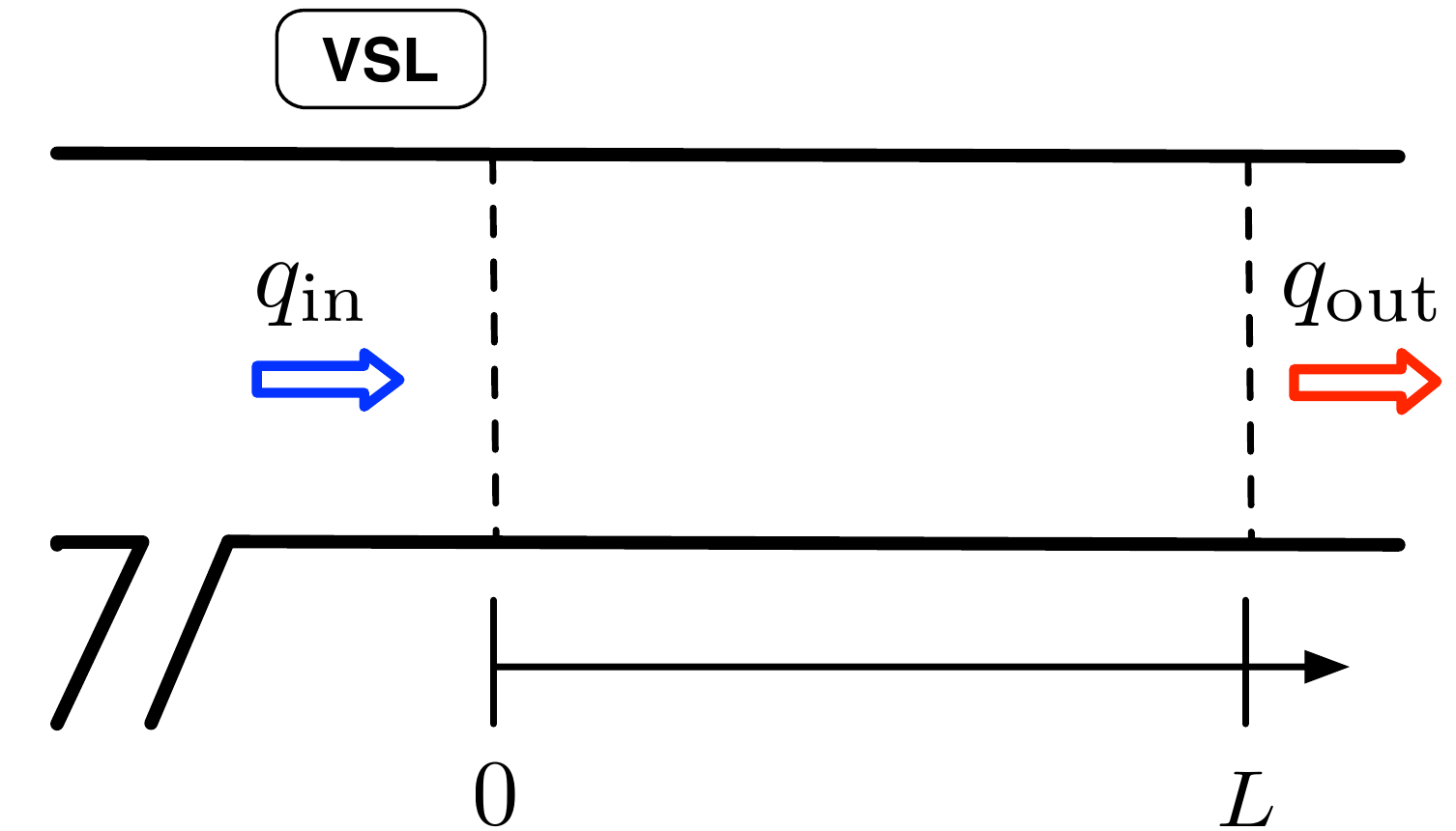}
	\caption{Traffic on a freeway segment}
\end{figure}

In practice, control of density at inlet can be realized with simultaneous operation of a ramp metering and a VSL at inlet as shown in Fig.3. The controlled density is then given by 
\begin{align}
\varrho(t) = \frac{q_{\textrm{in}}(t)}{v_c}.
\end{align}
where $v_c$ is the speed limit implemented by VSL and $q_{\textrm{in}}(t)$ is actuated by a on-ramp metering upstream of the inlet. Note that the linearized model is valid at the optimal density $\rho^{\star}$ since the reference density is assumed to be chosen near the optimal value.

\begin{figure*}[t]
	\includegraphics[width=\textwidth]{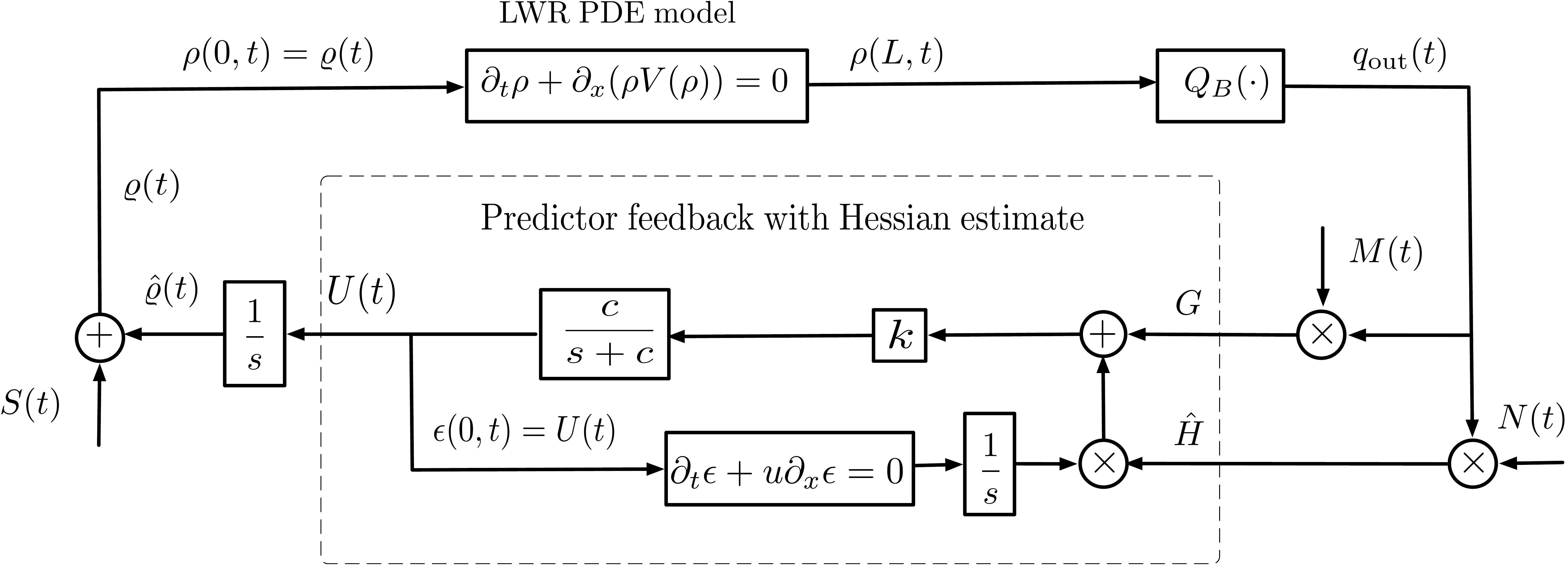}
	\caption{Block diagram for implementation of ES control design for nonlinear LWR PDE model}
	\label{fig:ES}
\end{figure*}
\section{Online Optimization by Extremum Seeking Control} 

In this section, we present the design of extremum seeking control with delay by following analogously the procedure in \cite{Tiago}. The block diagram of the delay-compensated ES algorithm applied to LWR PDE model is depicted in Fig. \ref{fig:ES}. 
Let $\hat{\varrho}(t)$ be the estimate  of $\rho^\star$, and $e(t)$ be the estimation error defined as
\begin{align}
e(t) = \hat \varrho(t) - \rho^\star. \label{e}
\end{align}
%
From Fig.4, the error dynamics can be written as
\begin{align}
		\dot e(t-D) = U(t-D). \label{doterror}
\end{align}
First, we introduce the dither signals $(M(t), N(t))$ given by 
\begin{align} 
M(t) =& \frac{2}{a} \sin\left( \omega t \right) , \\ N(t) =& - \frac{8}{a^2} \cos\left( 2 \omega t \right),
\end{align} 
where $a$ and $\omega$ are amplitude and frequency of a slow periodic perturbation signal $a \sin (\omega t)$ introduced later. 
Using the dither signals, we calculate estimates of the gradient and Hessian of the cost function, denoted as $(G(t), \hat{H}(t))$,
\begin{align} 
G(t) =& M(t) q_{\textrm{out}}(t) ,  \label{Gt} \\ \hat{H}(t) =& N(t) q_{\textrm{out}}(t),\label{Hhat} 
\end{align} 
where $\hat{H}(t)$ is to estimate the unknown Hessian $H$. The averaging of $G(t)$ and $\hat{H}(t)$ yields that 
\begin{align}
	G_{\rm av}(t) =& H e_{\rm av}(t-D), \label{Gav}\\
	\hat H_{\rm av} =& (Nq_{\textrm{out}})_{\rm av} = H. \label{HavH}
\end{align}
Taking average of \eqref{doterror}, we have
\begin{align}
		\dot e_{\rm av}(t-D) = U_{\rm av}(t-D),
\end{align}
where $ U_{\rm av}(t)$ is the averaged value for $U(t)$ designed later. Substituting the above equation into \eqref{Gav} gives that
\begin{align}
   \dot G_{\rm av}(t) = HU_{\rm av}(t-D). \label{tD}
\end{align}
The motivation for predictor feedback design is to compensate for the delay by feeding back future states in the equivalent averaged system $G_{av}(t+D)$. 
Given an arbitrary control gain $k>0$, we aim to design
\begin{align}
 U_{\rm av}(t) = k G_{\rm av}(t + D), \quad \forall t\geq 0. \label{tD1}
\end{align}
which requires knowledge of future states.
Therefore we have the following by plugging \eqref{tD1} into \eqref{doterror},
\begin{align}
\dot e_{\rm av}(t) =U_{\rm av}(t) = kH e_{\rm av}(t), \quad \forall t \geq D.
\end{align} 
Reminding that $k>0, H<0$, the equilibrium of the average system $e_{\rm av}(t) = 0$ is exponentially stable.

Applying the variation of constants formula $G_{\rm{av}}(t+D)=G_{\rm{av}}(t)+\hat{H}_{\rm{av}}(t) \int_{t-D}^{t}U_{\rm{av}}(\tau) d\tau$ and, from \eqref{tD1}, one has:
\begin{align}
	U_{\rm av}(t) =  k \left( G_{\rm av}(t) + \hat {H}_{\rm av}(t) \int_{t-D}^{t} U_{\rm av} (\tau) d \tau \right),  \label{Uavt}
\end{align}
which represents the future state $G_{\rm av}(t + D)$ in \eqref{tD} in terms of the average control signal $U_{\rm av}(\tau)$ for $\tau \in [t-D,t]$. The control input is infinite-dimensional due to its use of history over the past $D$ time units.

For the stability analysis in which the averaging theorem for infinite dimensional systems is used, we employ a low-pass filter for the above basic predictor feedback controller and then derive an infinite dimensional and averaging based predictor feedback given by
\begin{align} 
U(t) = {\mathscr  T} \left\{ k \left( G(t) + \hat{H}(t) \int_{t-D}^{t} U (\tau) d \tau \right) \right\} ,  \label{Ut}
\end{align} 
where $k>0$ is an arbitrary control gain, the Hessian estimate $\hat H(t)$ is updated according to \eqref{Hhat}, satisfying average property in \eqref{HavH}.  ${\mathscr T} \{\}$ is the low pass filter operator defined by 
\begin{align} 
{\mathscr  T} \left\{ \varphi (t) \right\} = {\mathscr L}^{-1} \left\{ \frac{c}{s+c} \right\} * \varphi(t), 
\end{align} 
where $c \in \mathbb{R^+}$ is the corner frequency, ${\mathscr L}^{-1} $ is the inverse Laplace transformation, and $*$ is the convolution in time.

\section{Stability Analysis} 
This section is devoted to the proof of the main theorem for delay-compensated ES algorithm, following \cite{Tiago}. Although the results in \cite{Tiago} are oriented to multiple and distinct delays, the derivation for the case of single delay was omitted there and it will be completely detailed here as a further contribution. 

 \begin{thm}
 	Consider the closed-loop system in Fig.4. There exits $c_0> 0$ such that $\forall c \geq c_0$, there exists $\omega_0(c_0)>0$ such that $\forall \omega > \omega_0$, the closed-loop system has a unique exponentially stable periodic solution in period $T = \frac{2\pi}{\omega}$, denoted by $e^T(t-D), U^T(\tau), \forall \tau \in [t-D,t]$, satisfying $\forall t> 0$
 \begin{align}
 		\left( |e^T(t-D)|^2 + |U^T(t)|^2 + \int_{0}^{D}|U^T(\tau)|^2 d\tau\right)^{\frac{1}{2}} \leq \mathcal{O}(1/\omega). \label{omega}
 \end{align}
 Furthermore,
\begin{align} 
	\lim\limits_{t \to +\infty} \sup|\varrho(t) - \rho^\star| =& \mathcal{O}(a+1/\omega),\\
	\lim\limits_{t \to +\infty} \sup|q_{\rm out}(t) - q^\star| =& \mathcal{O}(a^2+1/\omega^2). \label{limsup2}
\end{align}
\end{thm}
The proof of Theorem 1 is carried out along of the following Sections A to F.

\subsection{Closed-loop System} 
The estimate  $\hat{\varrho}(t)$ of the unknown optimal outgoing $\rho^\star$ is an integrator of the predictor-based feedback signal $U(t)$ as 
\begin{align} 
\dot{\hat{\varrho} }(t) = U(t), 
\end{align} 
and it follows that
\begin{align}
\dot 	e(t-D) = U(t-D).
\end{align}
The input $\varrho(t)$ to LWR PDE model is given by 
\begin{align} 
\varrho(t) = \hat{\varrho}(t) + S(t), \label{St}
\end{align} 
where the dither signal $S(t)$ is the inverse operator of a delayed perturbation signal $a \sin(\omega t)$, described as 
\begin{align} 
S(t) = a \sin \left( \omega (t+D) \right) .  \label{s}
\end{align} 
Substituting $S(t)$ into \eqref{St}, we have
\begin{align} 
\varrho(t) = \hat{\varrho}(t) + a \sin \left( \omega (t+D) \right).
\end{align} 
The delayed estimation error dynamics can be written as transport PDE system, $x \in [0,L]$
\begin{align} 
	\dot e(t-D) =& \epsilon(L,t), \label{PDE1}\\
\partial_{t} \epsilon(x,t) =& -  u \partial_{x} \epsilon(x,t)   , \\
{\epsilon}(0,t)= & U(t) . \label{PDEU}
\end{align} 
where
\begin{align} 
\epsilon(x,t) = U\left( t - \frac{x}{u} \right).  \label{solution_nonaverage}
\end{align} 
Combining \eqref{out}, \eqref{e}, and \eqref{s}, the relation among the estimation error $e(t)$, the input density $\varrho(t)$, and optimal outlet density $\rho^\star$ is given by
\begin{align}
	e(t) + a\sin(\omega t) = \varrho(t) - \rho^\star, \label{relation}
\end{align}
Substituting the above relation into the output map in \eqref{out}, we obtain the following equation
\begin{align}
q_{\textrm{out}}(t) 
= &  q^\star + \frac{H}{2} \left( e\left( t - D \right) + a \sin(\omega t) \right)^2. \label{qout}
\end{align}	
Plugging $M(t)$ and $G(t)$ into \eqref{Gt} and \eqref{Hhat} and representing the delayed input with PDE state $\epsilon(x,t)$, we have
\begin{align} 
U(t) =& {\mathscr  T} \left\{ k \left( G(t) + \hat{H}(t) \int_{0}^{L} \epsilon(\tau,t) d \tau \right) \right\} , \\
G(t) =& \frac{2}{a} \sin\left( \omega t \right)  q_{\textrm{out}}(t) ,   \\ \hat{H}(t) =&  - \frac{8}{a^2} \cos\left( 2 \omega t \right) q_{\textrm{out}}(t).
\end{align} 
It yields
{\small
\begin{align} 
U(t) =& {\mathscr  T} \left\{ k q_{\textrm{out}}(t)  \left( \frac{2}{a} \sin\left( \omega t \right) - \frac{8}{a^2} \cos\left( 2 \omega t \right)\int_{0}^{L} \epsilon(\tau,t) d \tau \right) \right\},
\end{align} }
and by substituting $q_{\textrm{out}}$ with \eqref{qout} and combining with transport PDE in \eqref{PDE1}-\eqref{PDEU}, we can write the closed-loop system as
\begin{align}
	\dot e(t-D) =& \epsilon(L,t), \label{close1}\\
	\partial_{t} \epsilon(x,t) =& -  u \partial_{x} \epsilon(x,t)   , \\
\notag 	{\epsilon}(0,t)= & {\mathscr  T} \Bigg\{ k \Bigg( q^\star + \frac{H}{2} \left( e\left( t - D \right) + a \sin(\omega t) \right)^2\Bigg) \\&  \Bigg( \frac{2}{a} \sin\left( \omega t \right) - \frac{8}{a^2} \cos\left( 2 \omega t \right)\int_{0}^{L} \epsilon(\tau,t) d \tau  \Bigg)  \Bigg\}.\label{close2}
\end{align}

\subsection{Average System}
Expanding \eqref{close2} and taking average of the closed-loop system, we obtain average model by setting the averages of sine and cosine functions of $n\omega, (n = 1,2,3,4)$ to zeros.
Note that the averaged controller satisfies 
\begin{align} 
\dot{U}_{av}(t) + c U_{av}(t) = c k \left( G_{av}(t) + H\int_{0}^{L} \epsilon_{\rm av} (\tau,t) d \tau \right),
\end{align} 
where $c>0$ is the corner frequency of the low pass filter and $k>0$ is the control gain. Denoting
\begin{align}
 \theta(t) = e(t-D), \label{PDE0_old}
\end{align}
the average system of \eqref{close1}-\eqref{close2} is rewritten by
\begin{align}
	\dot{{\theta}}_{\rm av}(t) =&
	\epsilon_{\rm av}(L,t), \label{av0}\\
	\partial_t \epsilon_{\rm av}(x,t) =&-u\partial_x \epsilon_{\rm av}(x,t), \label{av1}\\
 \partial_t\epsilon_{\rm av}(0,t) \!=&\! -c \epsilon_{\rm av}(0,t) \!+\! ckH \left(\theta_{\rm av}(t)\!+\!\int_{0}^{L}\!\!\epsilon_{\rm av}(\tau,t) d\tau \right). \label{av2}
\end{align}

\subsection{Backstepping Transformation} 
We apply backstepping transformation for the averaged delay state
\begin{equation}\label{Back_transf}
w(x,t) = \epsilon_{\rm av}(x,t) - kH\left[\theta_{\rm av}(t) + \int_{x}^{L}\epsilon_{\rm av}(\tau,t) d\tau\right].
\end{equation}
where $k>0$ and $H<0$. The average system \eqref{av0}-\eqref{av2} is mapped into the target system: 
\begin{eqnarray}
\dot{\theta}_{\rm av}(t)&=& kH\theta_{\rm av}(t)+ w(L,t), \label{targetsystem1}\\
\partial_t  w(x,t)&=& -u\partial_x  w(x,t) , \label{targetsystem2} \\
\partial_t w(0,t) &=& -(c + kH) w(0,t) \nonumber\\
&& - (kH)^2 \left[e^{\frac{kHL}{u}}{\theta}_{\rm av}(t) \right. \notag \\
&& \left. + \int_{0}^{L} e^{\frac{kH(L-\tau)}{u}}w(\tau,t) d\tau\right]~ \,. \label{targetsystem3}
\end{eqnarray}
We explain how to derive \eqref{targetsystem3} in detail. Combining \eqref{av2} and \eqref{Back_transf}, we have
\begin{align}
	w(0,t) =& -\frac{1}{c} \partial_t \epsilon_{\rm av}(0,t).
\end{align}
Taking time derivative on \eqref{Back_transf} for $w(0,t)$, we obtain 
\begin{align}
	\partial_t w(0,t) =& \partial_t \epsilon_{\rm av}(0,t) - kH \epsilon_{\rm av}(0,t).\label{wt}
\end{align}
The inverse transformation is given by
\begin{align}
\notag \epsilon_{\rm av}(x,t) =& w(x,t)+ kH\Bigg[e^{\frac{kH(L-x)}{u}}\theta_{\rm av}(t) \\ &+ \int_{x}^{L}e^{\frac{kH(L-x+\tau)}{u}}\epsilon_{\rm av} (\tau,t) d\tau\Bigg]. \label{inv}
\end{align}
Plugging \eqref{inv} and \eqref{av2} into \eqref{wt}, we obtain \eqref{targetsystem3} in the target system.
\subsection{Lyapunov Functional}
Now consider the following Lyapunov functional for the target system
\begin{equation}\label{Lyapunov}
V(t)\!=\!\frac{a{\theta}^2_{\rm av}(t)}{2}+\int_{0}^{L}e^{-x}w^2(x,t) dx + \frac{1}{2}w^2(0,t)\,,
\end{equation}
where the parameter $a>0$ is chosen later. Taking time derivative of the Lyapunov function, we have
\begin{align}
\notag	\dot{V}(t) =& akH \theta_{\rm av}^2 + a\theta_{\rm av} w(L,t)  + \frac{u}{2}w^2(0,t) -  \frac{ue^{-L}}{2}w^2(L,t)\\
\notag &-\frac{u}{2}\int_{0}^{L} e^{-x}w^2(x,t) dx + w(0,t) w_t(0,t)\\
\notag  \leq& akH \theta_{\rm av}^2 + \frac{a}{2b}\theta_{\rm av}^2 + \left(\frac{ab-ue^{-L}}{2} \right) w^2(L,t)\\
 &-\frac{u}{2}\int_{0}^{L} e^{-x}w^2(x,t) dx \notag\\
 &+ w(0,t) \left(w_t(0,t) + \frac{u}{2}w(0,t)\right) \label{vdot}
\end{align} 
where the positive constant $b$ satisfies the following,
\begin{align}
	b = \frac{ue^{-L}}{a}, \label{ab}
\end{align}
so that ${ab-ue^{-L}}  =0.$
The positive constant $a$ is chosen as 
\begin{align}
a = -ukHe^{-L}.\label{a}
\end{align}
Substituting $w_t(0,t)$ by \eqref{targetsystem3} and using Young's, Cauchy-Schwarz inequalities, the last term in \eqref{vdot} is bounded by
\begin{align}
\notag  &w(0,t) \left(w_t(0,t) + \frac{u}{2}w(0,t)\right)\\
\notag \leq& -\left(c - \frac{u}{2} + kH\right) w^2(0,t)\\ 
\notag & +\frac{e^L a^2}{4u}{\theta}^2_{\rm av}(t)  +  \frac{ue^{-L}}{a^2}\left|(kH)^2 e^{\frac{kHL}{u}}\right|^2 w(0,t)^2 \\ & + \frac{ue^{-L}}{4}\| w(t)\| ^2 + \frac{e^L}{u}\left \|(kH)^2 e^{\frac{kH(L-\tau)}{u}}\right \|^2 w(0,t)^2 
\label{lt}
\end{align}
Plugging \eqref{ab}--\eqref{lt} into \eqref{vdot}, one can arrive at 
\begin{align} \label{Lyapunov_derivative4}
\notag \dot{V}(t)\leq& - \frac{e^L a^2}{4u}{\theta}^2_{\rm av}(t) -\frac{u e^{-L}}{4}\int_{0}^{L} w^2(x,t) dx\\
& -(c - c_0) w^2(0,t),
\end{align}
where $c_0$ is defined as
\begin{equation} \label{cstar}
c_0=\frac{u}{2} - kH +  \frac{ue^{-L}}{a^2}\left|(kH)^2 e^{\frac{kHL}{u}}\right|^2 + \frac{e^L}{u}\left \|(kH)^2 e^{\frac{kH(L-\tau)}{u}}\right \|^2\, \\
\end{equation}
where $\tau \in [0,L]$.
An upper bound for $c_0$ can be obtained from lower and upper bounds of the unknown Hessian $H$.
Therefore, by choosing $c$ such that $c>c^*$, we obtain
\begin{equation} \label{Lyapunov_derivative5}
\dot{V}(t)\leq-\mu V(t)\,,
\end{equation}
for some $\mu>0$. Thus, the closed-loop system is exponentially stable in the sense of the $L^2$ norm
\begin{equation} \label{L2estimate_new}
\left(|{\theta}_{\rm av}(t)|^2 + \int_{0}^{L}w^2(x,t) dx + w^2(0,t)\right)^{1/2}.
\end{equation}
By the invertibility of the transformation, we can see that there exist constants $\alpha_1$ and $\alpha_2$ such that the following inequality is obtained 
\begin{eqnarray}
\alpha_1 \Psi(t) \leq  V(t) \leq \alpha_2 \Psi(t)\,,
\end{eqnarray}
where $\Psi(t)\triangleq |{\theta}_{\rm av}(t)|^2+\int_{0}^{L}
\epsilon^2_{\rm{av}}(x,t)dx + \epsilon^2_{\rm{av}}(L,t)$, or equivalently,
\begin{equation}
\Psi(t)\triangleq |{\theta}_{\rm av}(t-D)|^2+\int_{t-D}^{t}
U^2_{\rm{av}}(\tau)d\tau +U^2_{\rm{av}}(t)\,. 
\end{equation}
%
Hence, with (\ref{Lyapunov_derivative5}), we get
\begin{eqnarray} \label{final_inequality}
\Psi(t)\leq \frac{\alpha_2}{\alpha_1} e^{-\mu t} \Psi(0),\label{es}
\end{eqnarray}
which completes the proof of exponential stability of the averaged system. 
 
\subsection{Averaging Theorem}
The closed-loop system is written as
\begin{align}
\dot e(t-D) =& U(t-D),\\
\dot U(t) =& -c U(t) + c \left\{ k \left( G(t) + \hat{H}(t) \int_{t-D}^{t} U (\tau) d \tau \right) \right\}.
\end{align}
Defining the state vector $z(t)$ as $z(t) = [e(t-D), U(t)]^T$, and noting that 
$\int_{t-D}^{t}U(\tau) d\tau= \int_{-D}^{0}U(t+\tau) d\tau$, 
we can write the dynamics of $z$ as a functional differential equation described by 
\begin{align}
   \dot	z(t) = f (\omega t, z_t),
\end{align}
where $z_t(\tau) = z(t+\tau)$ for $-D\leq\tau\leq 0$. According to \eqref{es}, the origin of the average closed-loop system with transport PDE is exponentially stable. Applying the averaging theorem for infinite dimensional systems developed in \cite{Hale}, for $\omega$ sufficiently large, \eqref{close1}-\eqref{close2} has a unique exponentially stable periodic solution around its equilibrium satisfying \eqref{omega}. 

\subsection{ Asymptotic Convergence to a Neighborhood of the Extremum $(\rho^\star, q^\star)$}

By using the change of variables (\ref{PDE0_old}) and then integrating both sides of (\ref{close1}) within the interval $[t, \sigma+D]$, we have:
%
\begin{equation} \label{changingvariablesold}
{\theta}(\sigma+D)={\theta}(t)+\int_{t}^{\sigma+D}\epsilon(L,s) ds \,.
\end{equation}
From (\ref{solution_nonaverage}), we can rewrite (\ref{changingvariablesold})  in terms of $U$, namely
\begin{equation} \label{changingvariablesold2}
{\theta}(\sigma+D)={\theta}(t)+\int_{t-D}^{\sigma}U(\tau) d\tau \,.
\end{equation}
We define 
\begin{eqnarray} 
\vartheta(\sigma)={\theta}(\sigma+D)\,, \quad \forall \sigma \in [t-D,t] \,.
\end{eqnarray}
Applying \eqref{changingvariablesold2} to the above equation, we get
%
\begin{equation} \label{changingvariablesUav}
{\vartheta}(\sigma)=\vartheta(t-D)+\int_{t-D}^{\sigma}U(\tau) d\tau\,, \quad \forall \sigma \in [t-D,t]\,.
\end{equation}
By applying the supremum norm in both sides of (\ref{changingvariablesUav}) and using Cauchy-Schwarz inequality, we have
\begin{small}
	\begin{align}\label{changingvariablesupernew}
	\sup_{t-D \leq \sigma \leq t}\left|\vartheta(\sigma)\right|=& \sup_{t-D \leq \sigma \leq t}\left|{\vartheta}(t-D) \right|+\sup_{t-D \leq \sigma \leq t}\left| \int_{t-D}^{\sigma}U(\tau) d\tau \right| \nonumber\\
	\leq& \sup_{t-D \leq \sigma \leq t}\left|{\vartheta}(t-D)\right|+\sup_{t-D \leq \sigma \leq t} ~\int_{t-D}^{t}\left|U(\tau)\right| d\tau  \nonumber\\
	\leq& \left|{\vartheta}(t-D)\right|+\int_{t-D}^{t}\left|U(\tau)\right| d\tau \nonumber\\
	\leq& \left|{\vartheta}(t-D)\right|+\left(\int_{t-D}^{t}\!\!\!\!\! d\tau \right)^{1/2}  \left(\int_{t-D}^{t}\!\!\!\!\left|U(\tau)\right|^2 d\tau \right)^{1/2} \nonumber\\
	\leq& \left|{\vartheta}(t-D)\right|+ \sqrt{D} \left(\int_{t-D}^{t} U^2(\tau) d\tau \right)^{1/2}.
	\end{align}
\end{small}
One can easily derive 
\begin{small}
	\begin{align} 
	\left|{\vartheta}(t-D)\right| \leq& \left(\left|{\vartheta}(t-D)\right|^{2}+ \int_{t-D}^{t} \!\!\!\!U^2(\tau) d\tau \right)^{1/2}\!\,,\label{agoravai1} \\
	\left(\int_{t-D}^{t} U^2(\tau) d\tau \right)^{1/2} \leq& \left(\left|{\vartheta}(t-D)\right|^{2}+ \int_{t-D}^{t} \!\!\!\!U^2(\tau) d\tau \right)^{1/2}\!\,. \label{agoravai2}
	\end{align}
\end{small}
$\!\!$By using (\ref{agoravai1}) and (\ref{agoravai2}), one has
\begin{small}
	\begin{align} 
\notag	&\left|{\vartheta}(t-D)\right|+ \sqrt{D} \left(\int_{t-D}^{t} \!\!\!\!\!U^2(\tau) d\tau \right)^{1/2} \\
\leq& (1+\sqrt{D})\left(\left|{\vartheta}(t-D)\right|^{2}
 + \int_{t-D}^{t}U^2(\tau) d\tau \right)^{1/2}.
	\end{align} 
\end{small}
$\!\!$From (\ref{changingvariablesupernew}), it is straightforward to conclude that
\begin{small}	
	\begin{align} \label{changingvariablefinal}	
\sup_{t-D \leq \sigma \leq t}\left|\vartheta(\sigma)\right| \leq (1+\sqrt{D})\left(\left|{\vartheta}(t-D)\right|^{2}
+ \int_{t-D}^{t}U^2(\tau) d\tau \right)^{1/2},
\end{align} 	
\end{small}
and thus
\begin{small}
	\begin{align}  \label{changingvariablefinal_lim1}
	\left|\vartheta(t)\right|
	\leq& (1+\sqrt{D})\left(\left|{\tilde\theta}(t-D)\right|^{2}+ \int_{t-D}^{t} \!\!\!\!\!U^2(\tau) d\tau \right)^{1/2}.
\end{align} 
\end{small}
The above inequality (\ref{changingvariablefinal_lim1}) can be given in terms of the periodic solution $\vartheta^{\Pi}(t-D)$, $U^{\Pi}(\sigma)$, $\forall \sigma \in [t-D,t]$ as follows
\begin{small}
	\begin{align}  \label{changingvariablefinal_lim2}
	\left|\vartheta(t)\right|
	\leq& (1+\sqrt{D})\left(\left|\vartheta(t-D) -\vartheta^{\Pi}(t-D) + \vartheta^{\Pi}(t-D)  \right|^{2} \right. \nonumber \\
	& \left. + \int_{t-D}^{t} \left[U(\tau) - U^{\Pi}(\tau) + U^{\Pi}(\tau)\right]^2 d\tau \right)^{1/2} \,.
	\end{align} 
\end{small}
Applying Young's inequality, the right-hand side of (\ref{changingvariablefinal_lim2}) and $\left|\vartheta(t)\right|$ can be majorized by
\begin{small}
	\begin{align} \label{changingvariablefinal_lim3}
	\left|\vartheta(t)\right|
	\leq& \sqrt{2}~(1+\sqrt{D})\left(\left|{\vartheta}(t-D) -\vartheta^{\Pi}(t-D)\right|^2 + \left|\vartheta^{\Pi}(t-D)  \right|^{2} \right. \nonumber \\
	& \left. + \int_{t-D}^{t} \!\left[U(\tau) - U^{\Pi}(\tau)\right]^2 d\tau + \int_{t-D}^{t} \!\!\left[U^{\Pi}(\tau)\right]^2 d\tau\right)^{1/2}.
	\end{align}
\end{small}
$\!\!$From the averaging theorem \cite{Hale}, we have the exponential convergence
\begin{align}
	{\vartheta}(t-D)\!-\!\vartheta^{\Pi}(t-D) &\to 0\\
	\int_{t-D}^{t} \left[U(\tau)\!-\!U^{\Pi}(\tau)\right]^2 d\tau & \to 0
\end{align} 
Hence,
\begin{small}
	\begin{align}\label{limsupnonavg}
	\limsup_{t\to+\infty}|\vartheta(t)| =& \sqrt{2}~(1+\sqrt{D}) \nonumber \\
	& \times \left(\left|\vartheta^{\Pi}(t-D)\right|^{2}+ \int_{t-D}^{t} [U^{\Pi}(\tau)]^2 d\tau \right)^{1/2} .
	\end{align}
\end{small}
From (\ref{omega}) and (\ref{limsupnonavg}), we can write
\begin{align}
	\limsup_{t\to+\infty}|\vartheta(t)|\!=\!\mathcal{O}(1/\omega).
\end{align} 
From (\ref{e}) and recalling that $\varrho (t)= \hat \rho(t)+a\sin(\omega (t+D))$ and $\theta(t) = e(t-D) $, one has that
\begin{align} \label{error}
 \varrho(t) - \rho^\star=\vartheta(t) + a\sin(\omega (t+D)) \,. 
\end{align}
Since 
the first term in the right-hand side of (\ref{error}) is ultimately of order $\mathcal{O}(1/\omega)$ 
and the second term is of order $\mathcal{O}(a)$, then
%
%
\begin{equation}\label{limsup1new}
\limsup_{t\to+\infty}| \varrho(t) - \rho^\star|=\mathcal{O}(a+1/\omega)\,.
\end{equation}
Finally, from (\ref{out}), we get (\ref{limsup2}) and the proof is complete. $\hfill \Box$


\section{Simulation Result}
In simulation, we choose Greenshield's model for equilibrium velocity-density relation. For clear section Zone C, the fundamental diagram of traffic flow-density relation is given by
\begin{align}
Q(\rho) = -\frac{v_f}{\rho_m}\rho^2 + v_f \rho.
\end{align}
The maximum density is chosen to be 
\begin{align}
	\rho_m = \frac{5 \;\rm lanes}{7.5 \;\rm m } = 0.8 \;\rm vehicles/m,
\end{align}
where the $7.5 \;\rm m$ equals to the average vehicle length $5 \;\rm m$ plus $50 \% $ safety distance.
The maximum velocity is $v_f = 40 \rm \; m/s = 144 \; \rm km/h$. This $Q(\rho)$ is used in the nonlinear LWR PDE model simulation which describes the traffic dynamics upstream of bottleneck area. The maximum output flow also known as road capacity of Zone C is  
\begin{align}
q_c = \max_{0\leq \rho\leq\rho_m} Q(\rho) = 8 \;\rm vehicles/s.
\end{align}
The fundamental diagram in bottleneck area $Q_B(\rho)$, optimal/critical density $\rho^\star$ and maximum output flow $q^\star$ are unknown in practical implementation. The following function and parameters are chosen only for simulation purpose.
For bottleneck section Zone B, we consider the situation that only $3$ out of $5$ lanes still function. As a result, the road capacity reduces and we define the capacity reduction rate as $C_d = 60 \%$ compared with Zone C. Thus the following fundamental diagram is considered
\begin{align}
Q_B(\rho) = C_d Q(\rho)= -\frac{v_f}{\varrho_m}\rho^2 + v_f \rho,
\end{align}
where $\varrho_m = 0.48 \;\rm vehicles/m$ is the maximum density for reduced lanes in the bottleneck area and the same maximum velocity $v_f = 40\;\rm m/s = 144 \;\rm km/h$ is considered. The length of freeway segment is $L = 100 \;m$. If we consider a linearized LWR for Zone C, the characteristic speed is
\begin{align} 
u =& Q^\prime ( \rho)|_{\rho =\rho_r } = 20 \;s.
\end{align} 
The time delay for input is $D = \frac{L}{u}=5 \;s$. The outgoing flow $q_{\textrm{out}}(t)$ of the bottleneck area is 
\begin{align}
\notag		q_{\textrm{out}}(t) =& Q_B(\rho(L,t))\\
=&q^\star + \frac{H}{2} \left( \varrho\left( t - D \right)  - \rho^{\star} \right)^2,
\end{align}
where the optimal/critical density $\rho^\star$ and maximum output flow $q^\star$ are
\begin{align}
\rho^\star &= \frac{1}{2} \varrho_m = 0.24 \;\rm vehicles/m,\\
q^\star &= C_d q_c = 4.8 \;\rm vehicles/s.
\end{align}
The Hessian is obtained by taking second derivative of $Q_B(\rho)$
\begin{align}
H = -\frac{2 v_f}{\varrho_m} = - 166.7.
\end{align}
The Godnov scheme is employed for simulation of nonlinear LWR PDE model, which is commonly used in traffic flow application. The method is derived from the solution of local Riemann problems. The road segment is divided into spatial cell $\Delta x$ and the solution is advanced in time step $\Delta t$, which satisfy the following CFL condition
\begin{align}
	u_{\rm max}\frac{\Delta t}{\Delta x} < 1,
\end{align}
where $u_{\rm max}$ is the maximum characteristic speed. We choose the spatial cell $\Delta x = 0.05 \;\rm m$ sufficiently small so that numerical errors are negligibly small relative to the errors of the model.
\begin{figure}[t!]
	\includegraphics[width=0.49\textwidth]{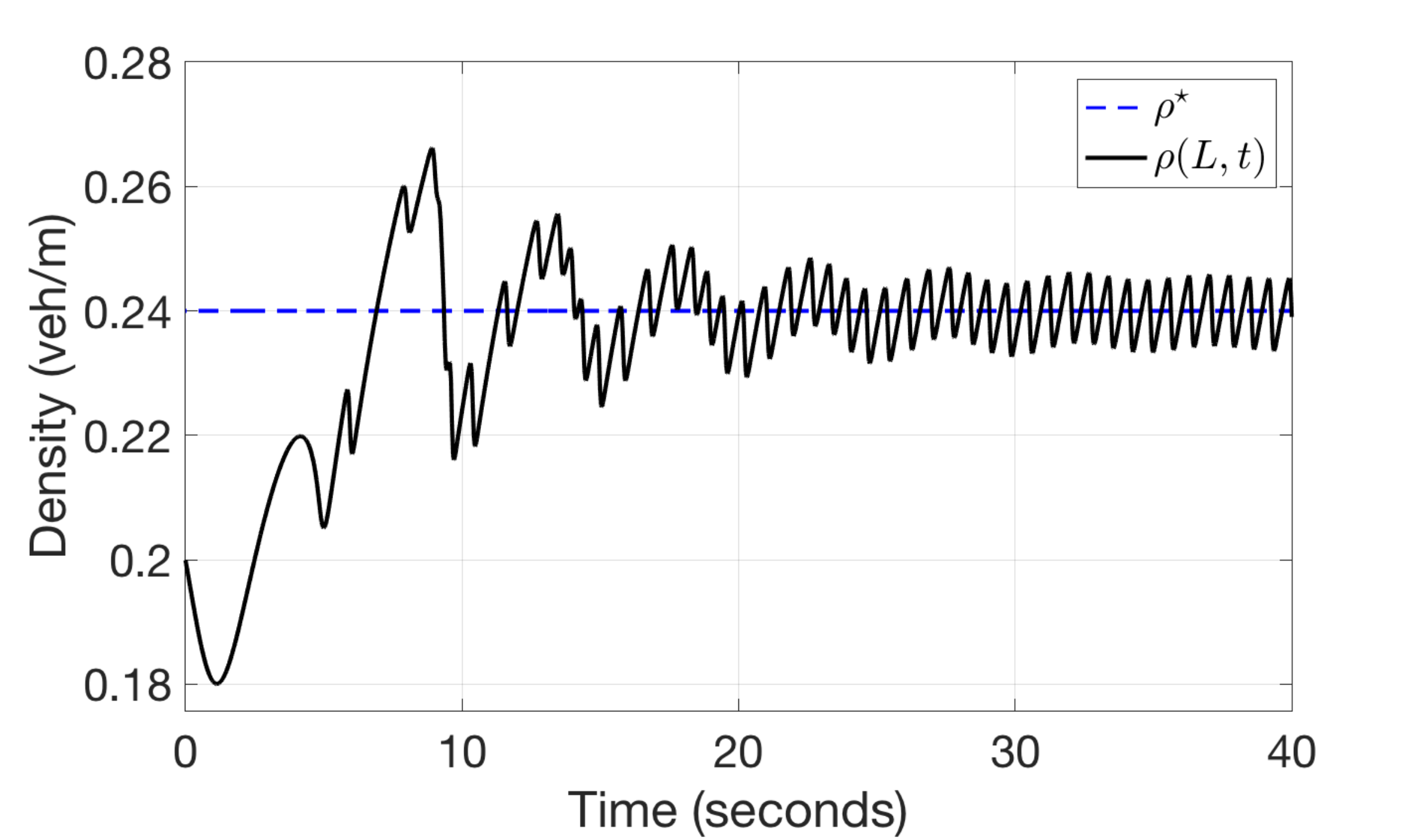}
	\caption{Traffic density $\rho(L,t)$ at the outlet of Zone C by nonlinear LWR model which is the input density for bottleneck area}
\end{figure}
\begin{figure}[t!]
	\centering
	\includegraphics[width=0.49\textwidth]{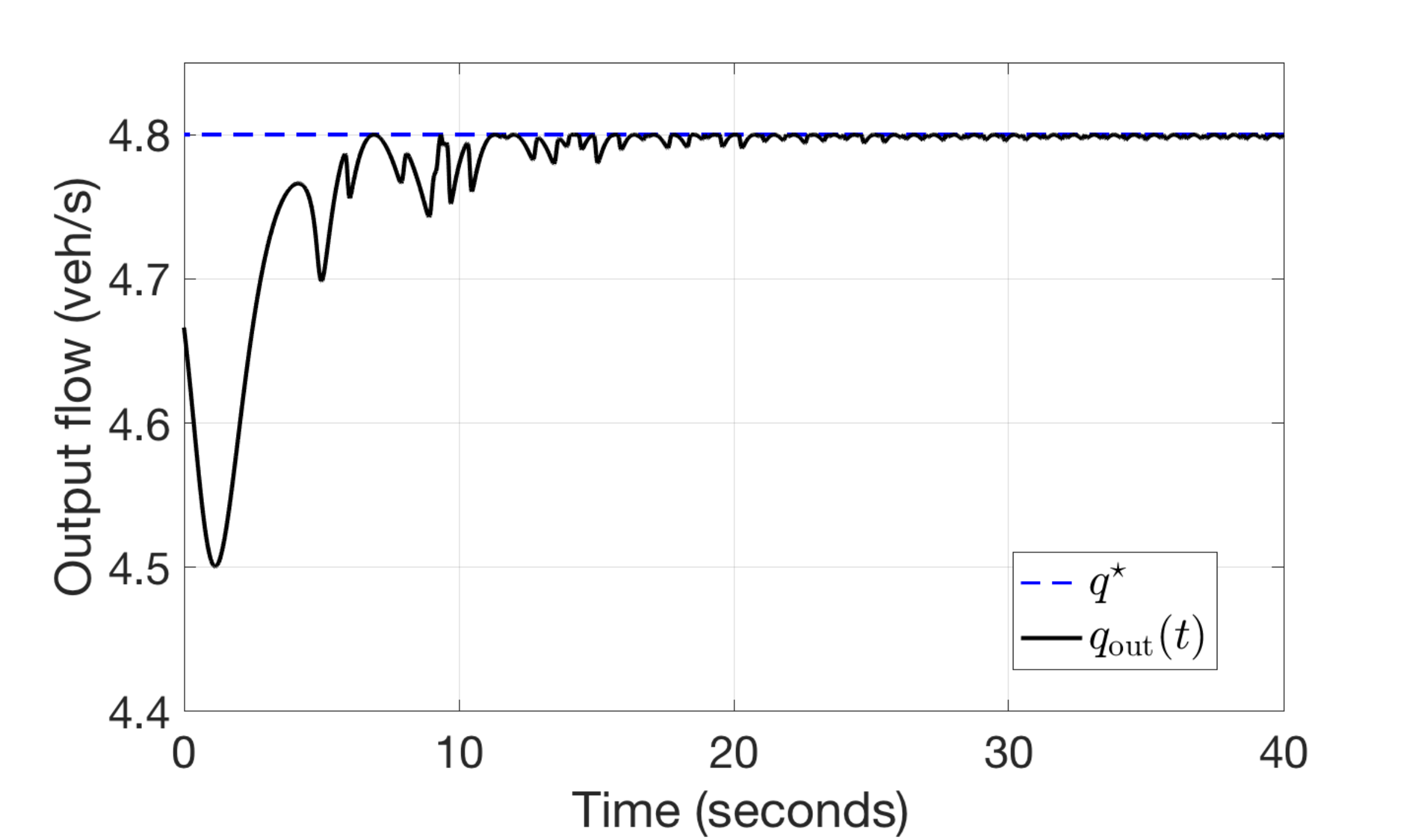}
	\caption{Outgoing traffic flow of the bottleneck area $q_{\rm out}(t)$ which is also the output flow for bottleneck area and the optimal value of outgoing flow $q^\star$}
\end{figure}
\begin{figure}[t!]
	\centering
	\includegraphics[width=0.49\textwidth]{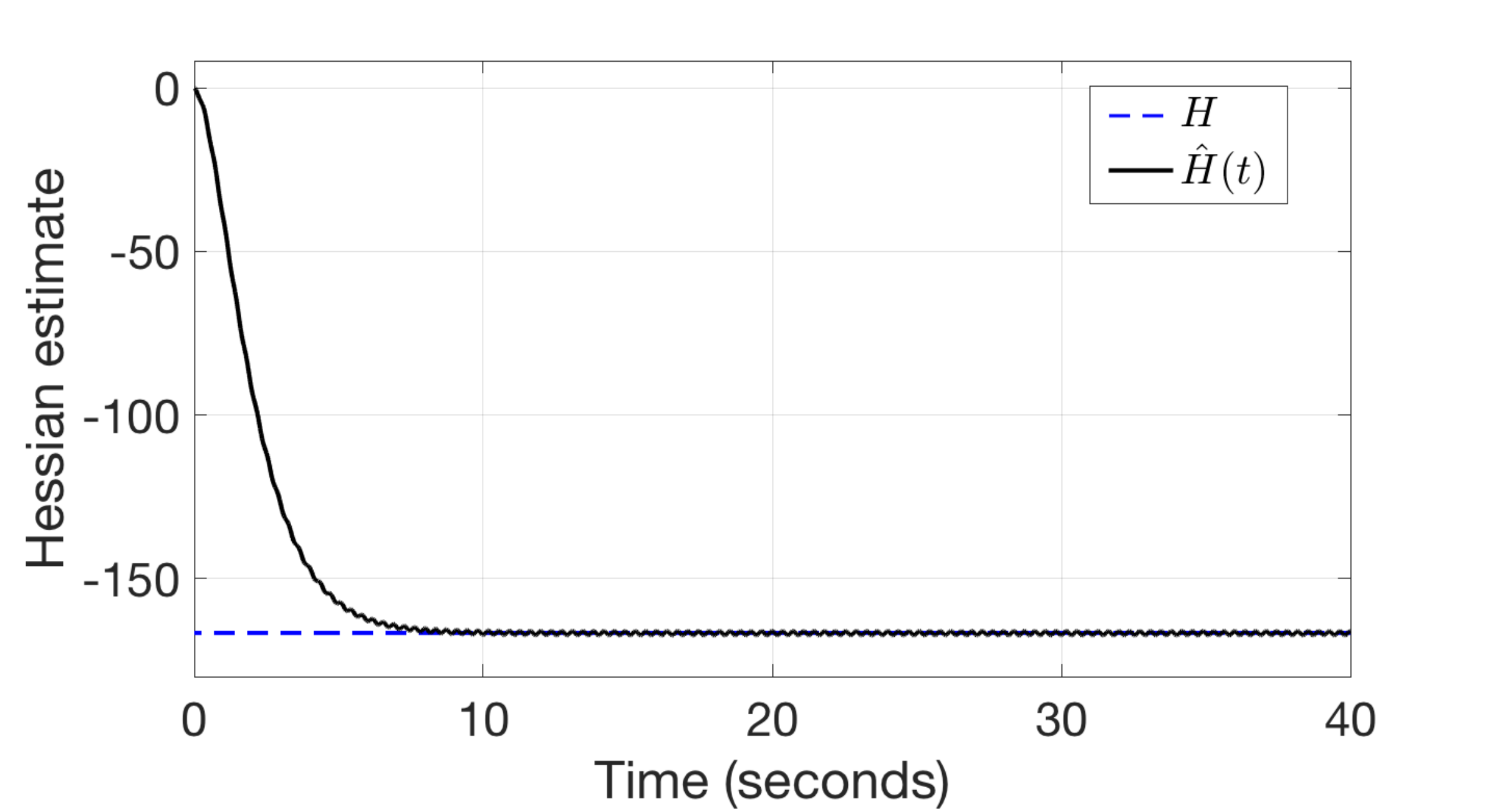}
	\caption{Hessian estimate $\hat H(t)$ of the ES control and prescribed Hessian value $H$}
\end{figure}

\begin{figure}[t!]
	\centering
	\includegraphics[width=0.51\textwidth]{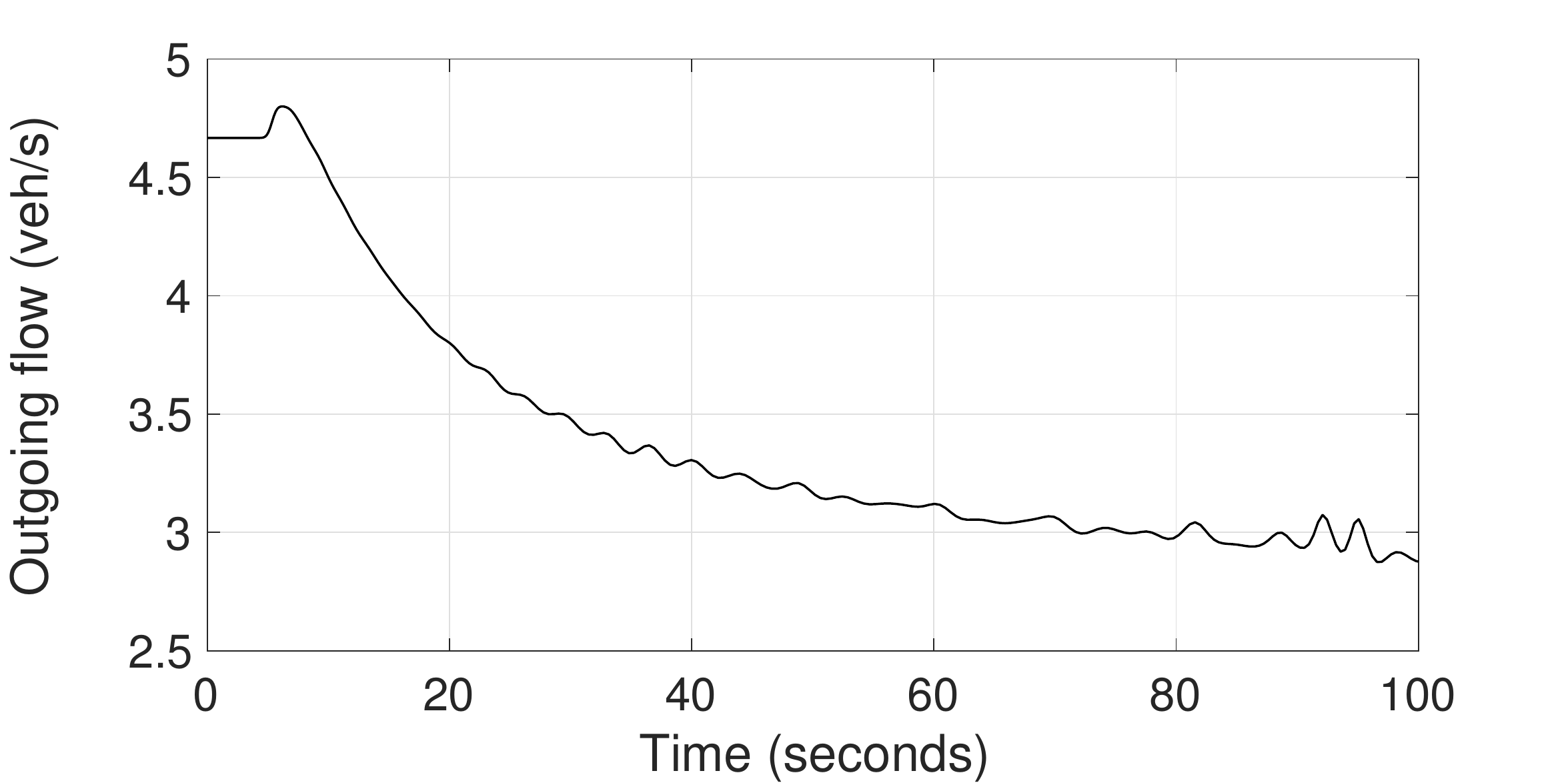}
	\caption{Output traffic flow of the bottleneck area without ES Control }
\end{figure}

The simulation result of the closed-loop system with ES control is shown in Fig.5-7. The parameters of the sinusoidal input and the designed controller are chosen to be $\omega = 2.75 \pi, a = 0.05, c = 50, K = 0.005$. One can observe that density in Fig.5 converges to a neighborhood of the optimal value $\rho^\star=0.24 \rm\; veh/m$ and the output flow of the bottleneck in Fig.6 converges to a neighborhood of the extremum point $q^\star = 4.8 \rm\; veh/s$. The Hessian estimate converges to the prescribed value $-166.7$. The convergence to optimal values is achieved in $40 \;\rm s$. In contrast, if we do not employ ES control for input density and the incoming flow depends only on upstream traffic. The open-loop system is shown in Fig.8. The evolution of outgoing flow at the bottleneck area is run for $100 \;\rm s$. We can see that the outgoing flow of the bottleneck area keeps decreasing and therefore congestion at the bottleneck area is getting worse till a bumper-to-bumper jam.

\section{Conclusion}
In this paper, we employ ES control to find a optimal density input for freeway traffic when there is a downstream bottleneck. To prevent traffic flow in bottleneck area overflowing the road capacity and furthermore causing congestion upstream in the freeway segment, the incoming traffic density at inlet of the freeway segment is regulated. The control design is achieved with delay compensation for ES control considering the upstream traffic is governed by the linearized LWR model. The optimal density and flow are achieved in the bottleneck area. The theoretical result is validated in simulation with the control design being applied on a nonlinear LWR PDE model along with an unknown fundamental diagram. Our future interest lies in conducting experimental validation of this problem. In a more sophisticated situation when there is multiple distant delays in presence of multi-lanes are going to be considered. It would also be interesting for authors to develop ES control with bounded update rates~\cite{Alex} under input delays exhibited through the LWR model and to develop a stochastic version of the algorithm presented in the paper by applying the results from \cite{Liu} and \cite{Rusiti}.


\begin{thebibliography}{99}  

\bibitem{Ariyur} 
K. B. Ariyur, and M. Krstic, 
\newblock {\em Real-Time Optimization by Extremum-Seeking Control}.  
\newblock John Wiley \& Sons, 2003. 

\bibitem{AW}  A. Aw, and M. Rascle, 
\newblock  `` Resurrection of "second order" models of traffic flow," 
\newblock {\em SIAM journal on applied mathematics}, vol.60, no.3, pp.916-938, 2000.

	
\bibitem{BL1} 
N. Bekiaris-Liberis, and M. Krstic,
\newblock  ``Compensation of transport actuator dynamics with input-dependent moving controlled boundary," 
\newblock {\em IEEE Transactions on Automatic Control}, vol.63, no.11, pp.3889-3896, 2018.

\bibitem{BL2}
N. Bekiaris-Liberis, and M. Krstic, 
\newblock ``Compensation of actuator dynamics governed by quasilinear hyperbolic PDEs,"
\newblock  {\em Automatica}, vol.92, pp.29-40, 2018.

\bibitem{benosman} 
M. Benosman, 
\newblock {\em Learning-Based Adaptive Control: An Extremum Seeking Approach--Theory and Applications,} 
\newblock Butterworth-Heinemann, 2016.

    
\bibitem{Jan} 
J. Feiling, S. Koga, M. Krstic, and T.R. Oliveira, 
\newblock ``Gradient extremum seeking for static maps with actuation dynamics governed by diffusion PDEs,"
\newblock {\em Automatica}, vol.95, pp.197-206, 2018.
	
\bibitem{Gha} 
A. Ghaffari, M. Krstic, and D. Nesic,
\newblock ``Multivariable Newton-based extremum seeking," 
\newblock {\em Automatica}, vol.48, no.8, pp.1759-1767, 2012.

\bibitem{Gha1} 
A. Ghaffari, M. Krstic, and  S. Seshagiri,
\newblock ``Power optimization and control in wind energy conversion systems using extremum seeking,"
\newblock {\em IEEE Transactions on Control Systems Technology}, vol.22, no.5, pp.1684-1695, 2014.

\bibitem{Gha2} 
A. Ghaffari, M. Krstic, and S. Seshagiri, 
\newblock ``Power optimization for photovoltaic microconverters using multivariable newton-based extremum seeking,"
\newblock {\em IEEE Transactions on Control Systems Technology}, vol.22, no.6, pp.2141-2149, 2014.




\bibitem{Guay2} 
M. Guay, D. Dochain, and M. Perrier, 
\newblock ``Adaptive extremum seeking control of continuous stirred tank bioreactors with unknown growth kinetics," 
\newblock {\em Automatica}, vol.40, no.5, pp.881-888, 2004.

\bibitem{Guay1} 
M. Guay, and T. Zhang, 
\newblock ``Adaptive extremum seeking control of nonlinear dynamic systems with parametric uncertainties,"
\newblock {\em Automatica}, vol.39, no.7, pp.1283-1293, 2003.

\bibitem{Hale} 
J.K. Hale, and S.V. Lunel, 
\newblock ``Averaging in infinite dimensions,"
\newblock {\em The Journal of integral equations and applications}, pp.463-494, 1990.

\bibitem{Jin}
H.Y. Jin, and W.L. Jin, 
\newblock ``Control of a lane-drop bottleneck through variable speed limits,"
\newblock {\em Transportation Research Part C: Emerging Technologies,} vol.58, pp.568-584, 2015.

\bibitem{Kan}
Y. Kan, Y. Wang, M. Papageorgiou, and I. Papamichail,
 \newblock  ``Local ramp metering with distant downstream bottlenecks: A comparative study,"
 \newblock {\em Transportation Research Part C: Emerging Technologies,} vol.62, pp.149-170, 2016.	
 
\bibitem{Liu} 
 S.J. Liu, and M. Krstic, 
\newblock ``Stochastic averaging in continuous time and its applications to extremum seeking,"
\newblock {\em IEEE Transactions on Automatic Control,} vol.55, pp.2235-2250, 2010.

\bibitem{MK1}
M. Krstic, 
\newblock {\em Delay Compensation for Nonlinear, Adaptive, and PDE Systems},
\newblock  Birkhäuser Boston, 2009.


\bibitem{MK6} M. Krstic, 
\newblock ``Performance improvement and limitations in extremum seeking control,"
\newblock  {\em Systems \& Control Letters}, vol.39, no.5, pp. 313-326, 2000.

\bibitem{MK5} 
M. Krstic, and H-H. Wang, 
\newblock ``Stability of extremum seeking feedback for general nonlinear dynamic systems,"
\newblock  {\em Automatica,} vol.36, pp.595-601, 2000.

\bibitem{MK4}
M. Krstic, and A. Smyshlyaev, 
\newblock {\em Boundary Control of PDEs: A Course on Backstepping Designs},
 \newblock Siam, 2008.

\bibitem{MK2} 
S.J. Liu, and M. Krstic,
\newblock {\em Stochastic Averaging and Stochastic Extremum Seeking}, 
\newblock Springer Science \& Business Media, 2012.


\bibitem{Tiago} 
T.R. Oliveira, M. Krstic, and D. Tsubakino, 
\newblock ``Extremum seeking for static maps with delays,"
\newblock  {\em IEEE Transactions on Automatic Control}, vol.62, no.4, pp.1911-1926, 2017.

\bibitem{Papa1} 
M. Papageorgiou, H. Hadj-Salem, and J.M. Blosseville,
\newblock ``ALINEA: A local feedback control law for on-ramp metering,"
\newblock {\em Transportation Research Record}, 1320(1), 58-67, 1991.

\bibitem{Paz}
P. Paz, T. R. Oliveira, A. V. Pino, and A. P. Fontana,
\newblock ``Model-Free Neuromuscular Electrical Stimulation by Stochastic Extremum Seeking,'' 
\newblock {\em IEEE Transactions on Control Systems Technology}, already on IEEE Xplore,  https://doi.org/10.1109/TCST.2019.2892924, 2019.

\bibitem{Rusiti} 
D. Ru{\v{s}}iti, G. Evangelisti, T.R. Oliveira, M. Gerdts, and M. Krstic, 
\newblock ``Stochastic extremum seeking for dynamic maps with delays,"
\newblock {\em IEEE Control Systems Letters,} vol.3, no.1, pp.61-66, 2019.

\bibitem{Alex}
A. Scheinker, and M. Krstic, 2014. 
\newblock ``Extremum seeking with bounded update rates,"
\newblock {\em Systems \& Control Letters,} vol.63, pp.25-31, 2014.

\bibitem{Papa} 
E. Smaragdis, M. Papageorgiou, and E. Kosmatopoulos, 
\newblock ``A flow-maximizing adaptive local ramp metering strategy,"
\newblock  {\em Transportation Research Part B: Methodological}, vol.38, no.3, pp.251-270, 2004.




\bibitem{Tan} 
Y. Tan, D. Ne{\v{s}}i{\'c}, and I. Mareels,
\newblock ``On non-local stability properties of extremum seeking control,"
\newblock {\em Automatica}, vol.42, no.6, pp.889-903, 2006.

\bibitem{Tan2} 
Y. Tan, D. Ne{\v{s}}i{\'c}, I. Mareels, and A. Astolfi, 
\newblock ``On global extremum seeking in the presence of local extrema,"
\newblock {\em  Automatica}, vol.45, no.1, pp.245-251, 2009.

\bibitem{Wang}
H. H. Wang, S. Yeung, and M. Krstic, 
\newblock ``Experimental application of extremum seeking on an axial-flow compressor,"
\newblock  {\em IEEE Transactions on Control Systems Technology}, vol.8, no.2, pp.300-309, 2000.

\bibitem{Wang} 
Y. Wang, E.B. Kosmatopoulos, M. Papageorgiou, and I. Papamichail, 
\newblock ``Local ramp metering in the presence of a distant downstream bottleneck: Theoretical analysis and simulation study,"
\newblock {\em  IEEE Transactions on Intelligent Transportation Systems,} vol.15, no.5, pp.2024-2039, 2014.

\bibitem{Huan}
H. Yu, and M. Krstic, 
\newblock ``Traffic congestion control for Aw-Rascle-Zhang model,"
\newblock  {\em Automatica,} vol.100, pp.38-51, 2019.	

\bibitem{Zhang}  H. M. Zhang,   
\newblock ``A non-equilibrium traffic model devoid of gas-like behavior,"
\newblock {\em Transportation Research Part B: Methodological,} vol.36, no.3, pp.275-290, 2002.
\end{thebibliography}
\end{document}